\documentclass{amsart}
\usepackage[utf8]{inputenc}
\usepackage{amsfonts,latexsym,amssymb,amsmath,amscd,amsthm}
\usepackage{graphicx}
\usepackage{enumerate}
\newtheorem{teorema}{Theorem}
\newtheorem{propo}[teorema]{Proposition}

\newtheorem{coro}[teorema]{Corollary}

\theoremstyle{definition}
\newtheorem{defin}[teorema]{Definition}

\theoremstyle{remark}
\newtheorem{nota}[teorema]{Remark}

\newtheoremstyle{teoremacita}% name of the style to be used
{3pt}% measure of space to leave above the theorem. E.g.: 3pt
{3pt}% measure of space to leave below the theorem. E.g.: 3pt
{\itshape}% name of font to use in the body of the theorem
{}% measure of space to indent
{\bfseries}% name of head font
{}% punctuation between head and body
{ }% space after theorem head; " " = normal interword space
{\thmname{#1}\thmnumber{ #2'}\thmnote{ #3}.}
\theoremstyle{teoremacita} \newtheorem*{teor*}{}

\newcommand{\be}{\begin{enumerate}}
\newcommand{\ee}{\end{enumerate}}

\newcommand{\bi}{\begin{itemize}}
\newcommand{\ei}{\end{itemize}}

\newcommand{\lista}[2]{{#1}_1,{#1}_2,\ldots ,{#1}_{#2}} %% Para escribir una lista 1, 2, ... , m

\newcommand{\ZZ}{{\mathbb Z}}
\newcommand{\NN}{{\mathbb N}}
\newcommand{\QQ}{{\mathbb Q}}
\newcommand{\CC}{{\mathbb C}}
\newcommand{\FF}{{\mathcal F}}

\newcommand{\eps}{\varepsilon}

\newcommand{\bx}{\mathbf{x}}

\newcommand{\bo}{{\mathbf 0}}

\newcommand{\cen}[1]{(\CC^{#1},\bo)}

\DeclareMathOperator{\Rs}{Rs}

\renewcommand{\theenumi}{\arabic{enumi}}
\renewcommand{\labelenumi}{\theenumi.}

\begin{document}
\title{Dicritical nilpotent holomorphic foliations}

%Autor: YO
\author{Percy Fern\'{a}ndez-S\'{a}nchez}
\address[Percy Fern\'{a}ndez]{Dpto. Ciencias - Secci\'{o}n Matem\'{a}ticas, Pontificia Universidad Cat\'{o}lica del Per\'{u}, Av. Universitaria 1801,
San Miguel, Lima 32, Peru}
\email{pefernan@pucp.edu.pe}

\author{Jorge Mozo-Fern\'{a}ndez}
\address[Jorge Mozo-Fern\'{a}ndez]{Dpto. \'{A}lgebra, An\'{a}lisis Matem\'{a}tico, Geometr\'{\i}a y Topolog\'{\i}a \\
Facultad de Ciencias, Universidad de Valladolid \\
Campus Miguel Delibes\\
Paseo de Bel\'{e}n, 7\\
47011 Valladolid - Spain}
\email{jmozo@maf.uva.es}
%\thanks{Gracias Esp\'{\i}ritu Santo por los favores recibidos}

%Otro autor
\author{Hern\'{a}n Neciosup}
\address[Hern\'{a}n Neciosup]{Dpto. Ciencias - Secci\'{o}n Matem\'{a}ticas, Pontificia Universidad Cat\'{o}lica del Per\'{u}, Av. Universitaria 1801,
San Miguel, Lima 32, Peru}
\email{hneciosup@pucp.pe}

\thanks{This work was funded by the Dirección de Gestión de la
Investigación at the PUCP through grant DGI-2015-1-0045, and by the Ministerio de Econom\'{\i}a y Competitividad from Spain, under Projects ``\'{A}lgebra y Geometr\'{\i}a en Din\'{a}mica Real y Compleja III" (Ref.: MTM2013-46337-C2-1-P) and ``\'{A}lgebra y geometría en sistemas dinámicos y foliaciones singulares'' (Ref: MTM2016-77642-C2-1-P)}

%\subjclass[2000]{Primary}


\date{\today}

\dedicatory{Dedicated to Felipe Cano, on the occasion of his first $60$ years of mathematical life}

\maketitle

\begin{abstract}
We study in this paper several properties concerning singularities of foliations in $\cen{3}$ that are pull-back of dicritical foliations in $\cen{2}$. Particularly, we will investigate the existence of first integrals (holomorphic and meromorphic) and the dicriticalness of such a foliation. In the study of meromorphic first integrals we follow the same method used by R. Meziani and P. Sad in dimension two. While the foliations we study are pull-back of foliations in $\cen{2}$, the adaptations are not  straightforward.
\end{abstract}

\section{Introduction}

In the paper \cite{FM}, first and second authors studied the analytic
classification of a generic class of quasi-ordinary
singularities of codimension one, holomorphic foliations
in dimension three. More precisely, we considered foliations
${\mathcal F}$ that can be generated by a holomorphic 1-form
$$
\Omega =d( z^2+x^py^q) + (x^{p'}y^{q'})^k \alpha U(x^{p'}y^{q'}) dz,
$$
where $p$, $q\in \NN \setminus  \{ 0\}$, $d=\gcd (p,q)$, $p=dp'$,
$q=dq'$, $U(t)\in \CC \{ t\}$, $U(0)=1$, $\alpha\in \CC^\ast$. In that paper, we focused  in foliations of generalized surface type, condition that implies certain relations between $k$ and $d$. We supposed that one of the two following conditions is satisfied:
\be[(i)]
\item Either $2k>d$,
\item Or $2k=d$ and a certain arithmetical condition on $\alpha$ is satisfied.
\ee
The classifying object turns out
to be the projective holonomy  of a certain component of the
exceptional divisor that appears in  a particular reduction of
the singularities, that follows the idea of Jung's method \cite{Giraud,Cossart}. This reproduces the scheme established in
dimension two by D. Cerveau and R. Moussu \cite{CM}, and R. Meziani \cite{Meziani}.

We are now  interested in the dicritical case, and the problem of knowing if such a foliation has a meromorphic first integral. As we will see, in this situation we must have that $d=2k$. Under this last condition,
several possibilities may appear. For generic values of  $\alpha$
the analytic classification of such foliations is essentially the
same that in the case $d<2k$, in the sense that projective
holonomy classifies.    We will suppose that we are in a
non-generic case, and that the foliation is dicritical. The
objective of the paper is to study  the existence, in
this case, of first integrals, both of holomorphic or  meromorphic type, and the dicriticalness of such a foliation. In particular, in order to investigate the existence of meromorphic first integrals we will follow the ideas developed in dimension two by
R. Meziani and P. Sad in \cite{MS}. Let us observe that the case
we are studying is a pull-back of a two-dimensional foliation.
More precisely, consider the map
\begin{eqnarray*}
\rho :\CC^3 & \longrightarrow & \CC^2 \\
(x,y,z) & \longmapsto & (x^{p'}y^{q'}, z).
\end{eqnarray*}
Then,
$$
\Omega= \rho^\ast (d(z^2+t^{2k}) + t^k \alpha U(t) dz).
$$
It is tempting to try to use this map in order to transfer the
properties of $\mathcal{F}_{k,\alpha,U}$ (foliation in $\cen{2}$
generated by $\omega_{k,\alpha,U}=d(z^2+t^{2k}) + t^{k} \alpha U(t) dz$) to
properties of ${\mathcal F}$. In fact, if there is a meromorphic
first integral $F(t,z)$ for $\mathcal{F}_{k,\alpha,U}$ it is clear that
$\Omega= \rho^\ast \mathcal{F}_{k,\alpha,U}$ has $\rho^\ast F = F\circ
\rho$ as meromorphic first integral. But in general it is not
true that  the existence of a first integral for a pull-back
implies that the original foliation has also a first integral.
Consider, e.g., any foliation generated by a 1-form $\omega$ in
$(\CC^2,\bo)$, and the restriction of its reduction of singularities to a small neighbourhood of a regular point. This restriction has always a first integral, while the original foliation may not. For the same reason, it is not straightforward that, if a foliation is dicritical, a pull-back of the foliation is also dicritical. The same example applies. 
Nevertheless, in the case under consideration, the condition for the dicriticalness and for the existence of
a meromorphic first integral is  essentially the same as in
dimension two. We will explore carefully this approach in this paper.

Let us briefly recall the criteria for the existence of a
meromorphic first integral for nilpotent foliations in
$(\CC^2,\bo )$. First of all, we have to remind that there is not
a topological characterization, as Suzuki example shows (see
\cite{CerveauMattei}). In our case, consider a foliation generated by
$$
\omega =d(z^2+t^{2p})+\alpha t^p U(t)dz,
$$
such that, after $p$ blow-ups, we arrive at a situation as
described in \cite{MS}: there are two singularities, not corners,
in the last component $\mathcal{D}_p$ of the divisor, call them $m_i$,
$m_{-i}$, such that  $m_i$ is imposed to be linearizable and dicritical, and
$m_{-i}$ simple resonant. Take a transversal $\Sigma$ on $\mathcal{D}_p$,
where the holonomy of $\mathcal{D}_p$ is computed. According to \cite{MS},
the existence of a meromorphic first integral  for $\omega$ is
equivalent to the existence of a rational function $r(T)$ on the
transversal $\Sigma$ invariant by the holonomy group of $\mathcal{D}_p$.
Besides this, Camacho \cite{Camacho} characterizes in full generality the existence of a
meromorphic first integral in terms of the singular holonomy of
the foliation. Let us briefly recall this notion and make the
link between both approaches, for the sake of completeness. The result shown by Camacho is the following one:
\begin{teor*}[Theorem \cite{Camacho}]
The foliation $\FF$ admits a meromorphic first integral if and only if:
\be[(i)]
\item The singular holonomy of the nondicritical branches is finite.
\item For each dicritical component $\mathcal{D}_{i}$ of the resolution of $\FF$, there is a meromorphic function $f_i:\mathcal{D}_{i}\rightarrow \overline{\CC}$ such that it is invariant by the singular holonomy of the nondicritical branches intersecting $\mathcal{D}_{i}$.
\item Given any two dicritical components $\mathcal{D}_i$ and $\mathcal{D}_j$, and a branch $B_{ij}$ between them, i.e., $q_i=B_{ij}\cap \mathcal{D}_{i}\neq \emptyset\neq B_{ij}\cap \mathcal{D}_{j}=q_{j}$, the germs of $f_{i}$ at $\mathcal{D}_{i}$, $f_{j}$ at $\mathcal{D}_{j}$ are compatible, in the sense that they take the same values ove the same leaves $L$ near $B_{ij}$.
\ee
\end{teor*}

Consider a corner of the divisor obtained after reduction of
singularities, where the foliation is resonant linearizable,
locally defined by $\omega_{pq}=pydx+qxdy$, $(p,q)=1$, and where
$\mathcal{D}_i$ is the component $x=0$, and $\mathcal{D}_j$ is $y=0$. The holonomy
$h_i$ of $\mathcal{D}_i$ at this singularity is $h_i(t)=e^{-2\pi ip/q}t$,
and analogously, $h_j(t)=e^{-2\pi iq/p}t$ turns out to be the
holonomy map of $\mathcal{D}_j$. Dulac's map is, in these coordinates, a map $D:
\Sigma_i \rightarrow \Sigma_j$ such that if $(t,1)\in \Sigma_i$,
$(1,D(t))\in \Sigma_j$ ($\Sigma_i$, $\Sigma_j$ are transversals
through $t=1$ to $\mathcal{D}_i$, $\mathcal{D}_j$, respectively), and $(t,1)$, $(1, D(t))$ are on the same leaf. So, we have
$$
t^p=D(t)^q\ \text{and then } D(t)=t^{p/q}.
$$
Define the adjunction map by $D\circ h_i=h_i^D \circ D$, $h_i^D:
\mathcal{D}_i\rightarrow \mathcal{D}_j$. A simple computation shows that
$h_i^D(t)=\xi t$, where $\xi^q=1$. The \textit{singular holonomy} group is the maximal subgroup obtained by iteration of the adjunction process made at the corners of a branch. In our case, as
$e^{-2\pi i p/q}$ is a primitive $q$-root of unity,
$\xi=(e^{-2\pi ip/q})^m$, and so, $h_i^D=h_j^m$, for some $m\geq
1$. So, the adjunction of $h_i^D$ to the holonomy  group in this
case does not increase it.

Consider now the situation studied by \cite{MS}. After reduction
of singularities, removing the only dicritical component, there
are two chains of divisors and all the corners are resonant
linearizable. Take a transversal $\Sigma$ on $\mathcal{D}_p$, the component
where the separatrix lies. Iterating the above construction, we
have that the singular holonomy group agrees with the usual one
of $\mathcal{D}_p$. In particular, it is finite. If $P\in \Sigma$, the leaf
through $P$ cuts the dicritical component $\mathcal{D}$. Assume that
there exists a rational function $f$ on $\mathcal{D}$,
invariant by the singular holonomy, as in (ii) of the Theorem above. The previous construction
defines a value $r(P)$ on $P\in \Sigma$, and so, a rational
function $r(T)$ such that $r\circ h=r$, for $h$ an element of the
holonomy group of $\mathcal{D}_p$. This is the condition of \cite{MS}.

Conversely, if such a $r(T)$ exists, by previous considerations,
$r\circ h=r$ for every $h$ in the singular holonomy
group, and a rational function on $\mathcal{D}$ can be defined.
These considerations relate the conditions stated by Camacho and
the ones by Meziani and Sad for the existence of a meromorphic
first integral. We will try to reproduce the last one in our
three dimensional, quasi-ordinary case.

The plan of the paper is as follows. In Section \ref{singularidadessimples}, following \cite{CC,Cano} we will recall the main notions concerning pre-simple singularities of foliations in dimension three, focusing in the existence of a meromorphic first integral. In Sections \ref{integralesholomorfas} and \ref{integralesmeromorfas} we will investigate the existence of holomorphic and meromorphic first integrals, respectively, for our three-dimensional foliations, relating the existence of these objects with the existence of the corresponding objects in dimension two. Finally, in Section \ref{dicriticidad} we will study the dicriticalness of such a foliation.

\subsection*{Acknowledgments} The authors want to thank Felipe Cano Torres, from the Universidad de Valladolid, for his continuous support and encouragement, and for fruitful conversations during the preparation of this work. First and third author wants to thank  Universidad de Valladolid, and  second author wants to thank Pontificia Universidad Católica del Perú for the hospitality during several visits to both institutions.

\section{Remarks on pre-simple singularities}
\label{singularidadessimples}
In this section we will recall the notion of pre-simple and
simple singularities of holomorphic foliations, as defined in
\cite{CC} and \cite{Cano}. We will establish some results
concerning their rigidity and the existence of meromorphic first
integrals.

Let ${\mathcal F}$ be a germ of holomorphic foliation in
$(\CC^n,\bo)$, generated by an integrable 1-form $\omega$,
adapted to a normal crossings divisor $E$. Choosing a
coordinate system such that  $E$ is written as $\prod_{i\in A}
x_i=0$, a meromorphic generator of ${\mathcal F}$ is
$$
\omega= \sum_{i\in A} a_i\frac{dx_i}{x_i}+ \sum_{i\notin A} a_i
dx_i ,\ (\lista{a}{n})=1.
$$
The component $x_i=0$ of the divisor will be invariant by
${\mathcal F}$ if $x_i$ does not divide $a_i$. A component
$x_i=0$ is called \textit{dicritical} if it is not invariant.
Denote
$$
A^{\ast} =\{ i\in A;\ x_i\nmid a_i\}
$$
the set of indexes corresponding to non-dicritical components.

The \textit{dimensional type} of ${\mathcal F}$, denoted
$\tau$, is the minimum number of coordinates needed to write a
generator of ${\mathcal F}$. Equivalently, it is the codimension
of the vector space
$$
{\mathcal X} (\mathcal{F})= \{ X(0); \ X\in {\mathcal X}
(\CC^n,\bo),\ \omega (X)=0\} ,
$$
(where ${\mathcal X}
(\CC^n,\bo)$ denotes the space of germs of holomorphic vector fields at the origin in $\cen{2}$), according to the rectification theorem for regular vector fields.

The \textit{adapted order} of ${\mathcal F}$ at the origin is defined as
$$
\nu ({\mathcal F},E):= \min \{ \nu (a_i); \ 1\leq i\leq n\},
$$
where $\nu (a_i)$ denotes the (usual) order of $a_i$ as formal
power series. The \textit{adapted multiplicity} of ${\mathcal F}$ at $\bo$
is
$$
\mu ({\mathcal F},E)= \min \{ \nu (a_i);\ i\in A\} \cup \{ \nu
(a_i)+1; \ i\notin A\} .
$$

The \textit{resonance invariant} of $({\mathcal F},E)$, $Rs ({\mathcal
F},E)$ is:
\bi
\item $\Rs ({\mathcal F},E)=1$ if $\nu ({\mathcal F},E)= \mu
({\mathcal F},E)$ and ${\mathcal F}$ is \textit{radially
dicritical}, i.e., the new component of the exceptional divisor that appears after
blowing-up the origin is dicritical.
\item $\Rs ({\mathcal F},E)=2$ if $\nu ({\mathcal F},E)=\mu
({\mathcal F},E)$, $\Rs ({\mathcal F},E)\neq 1$ and, in some local
coordinates, there is a map $\Phi : A\rightarrow \ZZ_{>0}$ such that
$$
\nu \left( \sum_{i\in A} \Phi (i)a_i\right)  > \nu ({\mathcal
F},E).
$$
\item $\Rs ({\mathcal F},E)=0$ in other case.
\ei

\begin{defin}
The singularity at the origin is \textit{pre-simple} for $({\mathcal
F},E)$ if $E$ is non-dicritical and one of the following holds:
\be
\item Either $\nu ({\mathcal F},E)=0$.
\item Or $\nu ({\mathcal F},E)= \mu ({\mathcal F},E)=1$, $Rs
({\mathcal F},E)=0$, and the directrix, i.e., the vector space
defined by $\{ In_1(a_i)=0;\ i\in A\}$ (where $In_1$ denotes the first order terms) has codimension one and
normal crossings with $E$.
\ee
\end{defin}

According to \cite[Prop. 46]{Cano}, there are formal coordinates $(\lista{x}{n})$ such that, any pre-simple singularity of dimensional type $\tau$ is given by one of the following formal meromorphic  types:

\be
\renewcommand{\theenumi}{\textbf{\Alph{enumi}}}
\renewcommand{\labelenumi}{\theenumi.}
\item  $\omega=\sum_{i=1}^\tau \lambda_i \frac{dx_i}{x_i}$, with
$\prod_{i=1}^\tau \lambda_i \neq 0$.
\item $\omega= \sum_{i=1}^k p_i\frac{dx_i}{x_i}+\psi (x_1^{p_1}\cdots x_k^{p_k}) \sum_{i=2}^\tau
\alpha_i  \frac{dx_i}{x_i}$,
where $p_i$ are natural numbers, $p_1\neq 0$, $\psi (\bo )=0$, $\alpha_i\neq 0$ for $i=k+1,\ldots, \tau$. This will be called the
resonant case.
\item $dx_1 -x_1 \sum_{i=2}^k p_i
\frac{dx_i}{x_i}+x_2^{p_2}\cdots x_k^{p_k}\sum_{i=1}^\tau \alpha_i
\frac{dx_i}{x_i}$, with $(p_2,\ldots ,p_k)\in \NN^{k-1}\setminus \{ (0,\ldots, 0)\} $, $\alpha_i\in \CC$. This is the Dulac
case.
\ee

We will focus here in the three dimensional case: singularities of
holomorphic foliations are of dimensional type $\tau=2$ or 3. According to the previous classification, pre-simple singularities are:
\be
\item If $\tau=2$, pre-simple is equivalent to the non-nilpotence
of the linear part of the dual vector field.
\item If $\tau=3$, they correspond to one of the
formal meromorphic models A, B or C, as above.
\ee

%\be
%\item $\omega=\sum_{i=1}^3 \lambda_i \frac{dx_i}{x_i}$, with
%$\prod_{i=1}^3 \lambda_i \neq 0$.
%\item $\omega= \sum_{i=1}^3 p_i\frac{dx_i}{x_i}+\sum_{i=2}^3
%\psi_i (x_1^{p_1}x_2^{p_2}x_3^{p_3}) \frac{dx_i}{x_i}$,
%where $p_i$ are natural numbers, $p_1\neq 0$. This is called the
%resonant case.
%\item $dx_1 -x_1 \sum_{i=2}^3 p_i
%\frac{dx_i}{x_i}+x_2^{p_2}x_3^{p_3}\sum_{i=1}^3 \alpha_i
%\frac{dx_i}{x_i}$,
%with $(p_2,p_3)\in \NN^2\setminus \{ (0,0)\} $. This is the Dulac
%case.
%\ee
%\ee

According to \cite{CC,Cano}, any germ of foliation in dimension
three can be transformed into another one having only pre-simple
singularities, using an appropriate sequence of blow-ups with permissible centers. If the
foliation has a meromorphic first integral, all the pre-simple
singularities that appear after this process must also have
meromorphic first integral, and as a consequence, the leaves of the foliation must be
closed. This fact imposes some restrictions that we will explore
here.

\newcounter{casos}
\begin{list}{\bfseries Case \Alph{casos}: }{\setlength{\leftmargin}{1cm} \setlength{\itemsep}{1ex}\usecounter{casos}}
%\addtolength{\leftmargin}{5cm}
%\renewcommand{\theenumi}{\Alph{enumi}}
%\renewcommand{\labelenumi}{\textbf{Case \theenumi.}}
\item Consider a pre-simple singularity, formally
equivalent to $\sum_{i=1}^3 \lambda_i
\frac{d\hat{x}_i}{\hat{x}_i}$. It has a multiform first integral
$\hat{x}_1^{\lambda_1}\hat{x}_2^{\lambda_2}\hat{x}_3^{\lambda_3}$.
According to \cite[Cor. 1.2., pg. 115]{CerveauMattei} it is
semi-divergent and consequently the singularity has three
convergent separatrices, tangent to $x_i=0$. Consequently, in
appropriate analytic coordinates the foliation is defined by
$$
\sum_{i=1}^3 \lambda_i (1+b_i(\bx)) \frac{dx_i}{x_i},
$$
with $b_i (\bo)=0$. Consider the separatrix $x_1=0$. Its holonomy
group is generated by two diffeomorphisms $h_2(x_1)$, $h_3(x_1)$,
such that
$$
h_k (x_1)= e^{-\frac{\lambda_k}{\lambda_1} 2\pi i} x_1+\cdots .
$$
If the leaves of the foliation are closed, the holonomy is
periodical, i.e.,
$$
\frac{\lambda_k}{\lambda_1}=\frac{p_k}{p_1},\ p_i\in \ZZ.
$$
So, it is linearizable and then, the foliation also is, i.e. a
generator is
$$
\omega= \sum_{k=1}^3 p_k\frac{dx_k}{x_k},
$$
with $(p_1,p_2,p_3)\in (\ZZ\setminus \{ 0\})^3$. So, existence of
a meromorphic first integral for a pre-simple singularity of type
\textbf{A}, and dimensional type $\tau=3$ implies linearizability.

\item If $k=2$, the foliation is a saddle-node. Indeed, assume $p_3=0$, and take a transversal section $x_1=\varepsilon$, for $\eps\neq 0$ small enough. The resulting two-dimensional foliation is
$$
x_3 (p_2+\alpha_2 \psi (\eps^{p_1} x_2^{p_2})) dx_2 + \alpha_3 x_2 \psi (\eps^{p_2} x_2^{p_2}) dx_3,
$$
which is a saddle-node. Consequently, the leaves are not closed, and a meromorphic first integral cannot exist.

If $k=1$, and the \emph{residual spectrum of the singularity} \cite[Rem. 20]{Cano} is resonant, i.e., $\alpha_2/\alpha_3\in \QQ_{<0}$, after blowing-up the $x$-axis a certain number of times we arrive to a saddle-node again. If this residual spectrum is non resonant, or $k=3$, the singularity is simple.
In this case, the holonomy group of $x_1=0$ has been computed in \cite[Section 4]{CerveauMozo}. Under the assumption that the leaves are closed, this group turns out to be linearizable, and the foliation has a holomorphic first integral such that, in appropriate coordinates can be written as $x_1^{p_1}x_2^{p_2}x_3^{p_3}$, so, in fact, it is of the previous type \textbf{A}.

\item  An analytic generator of the foliation is
$$
\omega=dx_1-(p_2x_1+\cdots )\frac{dx_2}{x_2}-(p_3x_1+\cdots
)\frac{dx_3}{x_3}.
$$
Consider the following commuting  vector fields, tangent to
$\omega$:
\begin{align*}
X & = (p_2x_1+\cdots)\frac{\partial}{\partial x_1}+x_2
\frac{\partial}{\partial x_2} \\
Y & = (p_3x_1+\cdots)\frac{\partial}{\partial x_1}+x_3
\frac{\partial}{\partial x_3}.
\end{align*}
Take $b>0$ and define $X_b=X+bY$. The linear part of $X_b$ is
$$
X_{b,1}= (p_2+bp_3)x_1\frac{\partial}{\partial x_1}+
x_2\frac{\partial}{\partial x_2}+ bx_3 \frac{\partial}{\partial
x_3}.
$$
The triple $(q_1,q_2,q_3)\in \NN^3$ is a resonance for $X_b$ if
$$
q_1 (p_2+bp_3)+q_2+q_3 b=p_2 +bp_3.
$$
Assume $b\notin \QQ$. The previous equality means that
$$
q_1p_2+q_2=p_2;\qquad q_1p_3+q_3=p_3,
$$
and so, $q_1\in \{ 0,1\}$. If $q_1+q_2+q_3\geq 2$,
$(q_1,q_2,q_3)=(0,p_2,p_3)$, and by Poincar\'{e} normalization
theorem, $X_b$ is analytically equivalent to
$$
A_b= \left( ( p_2+bp_3)x_1+a_1 (x_2^{p_2}x_3^{p_3})\right)
\frac{\partial}{\partial x_1}+ x_2 \frac{\partial}{\partial x_2}
+ bx_3 \frac{\partial}{\partial x_3}.
$$
Consider $b'\neq b$, $b'\notin \QQ$, and construct analogously $A_{b'}$. The
1-form
$$
\tilde{\omega}=\iota_{A_b}\iota_{A_{b'}} (dx_1\wedge dx_2 \wedge
dx_3
)
$$
is a generator of the foliation and implies that the normal form
in this case is, in fact, convergent. $x_2=0$, $x_3=0$ are
separatrices. The holonomy generators $h_1(x_3)$, $h_2(x_3)$ of
$x_3=0$ verify
\begin{align*}
-\frac{1}{h_1(x_3)^{p_3}} +\alpha_3 \log
(h_2(x_3))+\alpha_2 & =-\frac{1}{x_3^{p_3}}+\alpha_3\log (x_3),
\\
\frac{1}{h_2(x_3)^{p_3}}- \alpha_3 \log(h_2(x_3))
 & =\frac{1}{x_3^{p_3}} -\alpha_3 \log(x_3).
\end{align*}

The leaves being closed, $\alpha_3=0$. This implies that
$$
{h_2(x_3)^{p_3}} ={x_3^{p_3}},\ \text{i.e.
}h_2(x_3)=\xi x_3,\ \xi^{p_3}=1.
$$
We also have that
$$
-\frac{1}{h_1(x_3)^{p_3}}+\alpha_2=- \frac{1}{x_3^{p_3}},
$$
and so,
$$
h_1(x_3)^{p_3}=\frac{x_3^{p_3}}{1+\alpha_2 x_3^{p_3}}.
$$
So, $h_{1}(x_3)=\eta  x_3+ \ast \cdot x_3^{p_3}+\cdots $, with
$\eta^{p_3}=1$. Periodicity of $h_1$ implies that $\alpha_2=0$.
This forbids singularities of type \textbf{C} with
$(\alpha_2,\alpha_3)\neq (0,0)$ to have a meromorphic first
integral.
\end{list}

\section{Holomorphic first integrals} \label{integralesholomorfas}

Let us study the problem of the existence of holomorphic first integrals. According to \cite{MM}, a singular holomorphic foliation has a first integral if and only if all the leaves are closed and only a finite number of them contain the singularity in their closure. So, it is a property of topological nature.

Let us consider a nilpotent singularity on $(\CC^{m+1}, \bo )$, i.e., a singularity of foliation of the type $zdz+\cdots $. Under a coordinate change, using Loray's Preparation Theorem \cite{Loray}, this foliation can be generated by a 1-form
$$
\Omega:=\Omega_{f,n,p,\alpha,U}=d(z^2+f^n) +\alpha f^p U(f)dz,
$$
where $f(\bx ) :(\CC^m,\bo )\rightarrow (\CC,0 )$, $\alpha\in \CC^\ast $, $U(t)\in \CC\{ t\}$, $U(\bo )=1$. Write $f=f_1^{r_1}\cdots f_k^{r_k}$, $\gcd (r_1,\ldots , r_k )=1$. Consider the map
$$
\begin{array}{rcl}
\Phi: (\CC^{m+1}, \bo ) & \longrightarrow & (\CC^2, \bo ) \\
(\bx, z) & \longmapsto & (f(\bx ), z) ,
\end{array}
$$
so, $\Omega=\Phi^\ast \omega_{n,p,\alpha,U}$, where
$$
\omega_{n,p,\alpha,U} =d(z^2+t^n) +\alpha t^p U(t) dz.
$$
It is clear that, if $h(t,z) $ is a holomorphic first integral
 for $\omega_{n,p,\alpha,U}$, $h\circ \Phi$ is a holomorphic first integral for $\Omega$. Let us show the converse in this case.

 The map $\Phi$ is a surjective submersion on $(\CC^{m+1},\bo )\setminus S$, where
 $$
 S:=\left\{ (\bx , z)\in (\CC^{m+1},\bo );\ \nabla f(\bx ) =\bo \right\} \subseteq S':= \left\{  (\bx , z)\in (\CC^{m+1},\bo );\ f(\bx )=0\right\} .
$$
\begin{propo}
If $\mathcal{L}$ is a closed leaf for $\Omega$, $\Phi (L) $ is a closed leaf for $\omega_{n,p,\alpha}$.
\end{propo}

\begin{proof}
$\mathcal{L}\setminus S'$ is closed on $(\CC^{m+1},\bo )\setminus S'$, and $ \Phi (\CC^{m+1} \setminus S')= \CC^\ast \times \CC$ (at least as a germ around the origin, we won't precise the neighbourhood in which this equality is valid). By \cite[Thm. 3, Ch. 2]{CLN}, $\mathcal{L}\setminus S'$ is a connected component of $\Phi^{-1} (F')$, $F'= \Phi (\mathcal{L}\setminus S' ) $ being a leaf of $\mathcal{F}_{\omega_{n,p,\alpha,U}}$ in $\CC^\ast \times \CC$.

Consider a trivializing chart $U$ for the foliation defined by $\omega_{n,p,\alpha,U}$ around $P\in F'$, where the foliation is defined by the levels of a submersion $g:U\rightarrow \CC$. For every point $(t_0,z_0)\in U$, the inverse image by $\Phi$ is the variety $f(\bx )=t_0$, $z=z_0$. This variety has a finite number of connected components, so, reducing $U$ if necessary, $\Phi^{-1} (U)$ consists in a finite number of connected trivializing sets. The leaf $\mathcal{L}\setminus S'$ being closed, in each of these open sets there are a finite number of plaques of $\mathcal{L}\setminus S'$, so $F$ has a finite number of components in $U$ and then, it is closed.

So, $\Phi (\mathcal{L}\setminus S')$ is closed on $ \CC^\ast \times \CC$. As
$$
\Omega \wedge df = (2z+ \alpha f^p U(f)) dz\wedge df,
$$
$S'$ is transversal to $\mathcal{F}_{\Omega}$ outside $z=0$. Consider now the leaf $F$, such that $\mathcal{L}$ is a connected component of $\Phi^{-1} (F)$. In fact, $F=\Phi (\mathcal{L})$. This leaf cuts the axis $t=0$ transversally outside $z=0$.
As $\mathcal{L}$ is closed, the ``vertical'' axis $z$ cuts transversally $\mathcal{L}$ in a discrete number of points, so $F$ does also. Then $F$ is closed outside the origin.
Apply finally Remmert-Stein Theorem in order to extend $F$ up to the origin, concluding that $\Phi (\mathcal{L}) $ is closed.

\begin{figure}
\includegraphics[width=10cm]{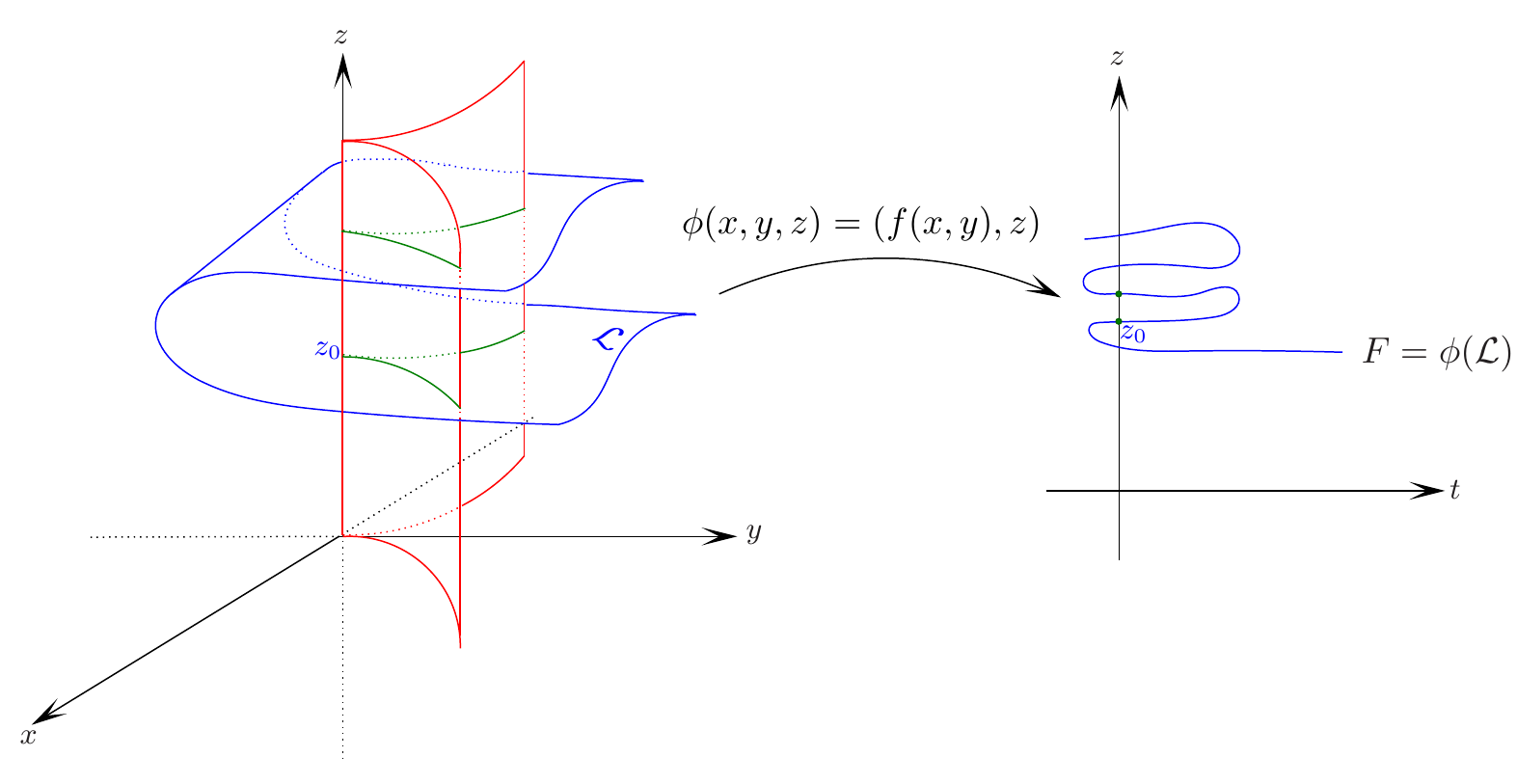}
\caption{The image of a closed leaf is closed}
\end{figure}

%So, $S\cap L$ has a finite number of connected components, and so $\Phi (S)\cap \{ t=0\}$ consists in a finite number of pints. Again by \cite[Ch. 3, Thm. 5]{CLN}, $\Phi (L)$ is closed on $\CC^2 \setminus \{ \bo\}$. Remmert-Stein Theorem allows to extend $\Phi (L)$ up to the origin.

\end{proof}

Consequently, if $\Omega$ has a holomorphic first integral, as $\Phi$ is surjective, the leaves of the foliation generated by $\omega_{n,p,\alpha,U}$ are closed. If there exists infinitely many separatrices, their inverse images would be separatrices of $\Omega$ containg the analytic set $\Phi^{-1}(\bo )$, which contradicts the fact that $\Omega$ has a holomorphic first integral. When a singular holomorphic foliation has a first integral, \cite{MM} allows to conclude the following result:
\begin{teorema}
The holomorphic foliation $\mathcal{F}_{\Omega_{f,n,p,\alpha,U}}$ in $(\CC^{m+1},\bo)$ has a holomorphic first integral if and only if $\mathcal{F}_{\omega_{n,p,\alpha,U}}$ has a holomorphic first integral.
\end{teorema}

In \cite{CM}, criteria are given for the existence of a holomorphic first integral in terms of the projective holonomy. The generators of this group must be of finite order, and also the group, which implies its abelianity. Combining these criteria with the results in \cite{Meziani}, we have:

\begin{coro}
A nilpotent singular foliation $\mathcal{F}_{\Omega}$ defined by $\Omega= \Omega_{f,n,p,\alpha,U}=0$ has a holomorphic first integral if and only if one of the following conditions occurs:
\be
\item $n<2p$, the projective holonomy of $\mathcal{F}_{\omega_{n,p,\alpha,U}}$ is abelian, and its generators have finite order.
\item $n=2p$, $\alpha^2\neq \dfrac{(16+r)^2}{16+2r}$, $\forall\, r\in \QQ_{\geq 0}$, and also the projective holonomy of $\mathcal{F}_{\omega_{n,p,\alpha,U}}$ is abelian, and its generators have finite order.
\ee
\end{coro}
Note that the conditions on $n$, $p$ are necessary to guarantee that no dicritical components and no saddle-nodes appear in the reduction of singularities of $\mathcal{F}_{\omega_{n,p,\alpha,U }}$.

Let us mention that in \cite{Moussu}, an example of a nilpotent foliation without a holomorphic first integral is presented. The projective holonomy of the foliation defined by the 1-form $\omega=d(y^2+x^3)+x (2xdt-3ydx)$
has two generators $h_1$, $h_2$, with $h_1^2=h_2^3=id$ but they don't commute, so the projective holonomy is not finite. 

\section{Meromorphic first integrals in the quasi-ordinary case}
\label{integralesmeromorfas}
Let us focus now in cuspidal, quasi-ordinary singularities of
holomorphic foliations in dimension three. A generator  in
appropriate coordinates is 
$$
\Omega=d(z^2+(x^py^q)^k)+\alpha (x^py^q)^n U(x^py^q)dz,
$$
where $\alpha\in \CC^\ast$, $(p,q)=1$, $U(t)\in \CC \{ t\}$,
$U(\bo )=1$ \cite{FM,FMN1}. In this  section we are interested in the existence of pure meromorphic first integrals (i.e., non holomorphic, as defined in \cite{CerveauMattei}). The following result restricts the values of $n$ and $k$ that we have to study:
\begin{propo}
If a quasi-ordinary singularity of foliation, generated by a 1-form
$$
d(z^2+(x^py^q)^k)+\alpha (x^py^q)^n U(x^py^q)dz,
$$
has a pure meromorphic first integral, then $k=2n$.
\end{propo}
\begin{proof}
If $2n>k$, it is a generalized surface. If $2n<k$, blow-up $pn$ times the $y$-axis, and $qn$ times the $x$-axis. The equations of this transform are,  in the main chart, $z=tx^{pn}y^{qn}$. The strict transform of $\omega$ is
$$
\tilde{\Omega}= (2nt^2 +k (x^py^q)^{k-2n} + \alpha nt U) (pydx+qxdy) +xy (2t +\alpha U) dt.
$$
The origin turns out to be a simple singularity of type \textbf{B}, resonant, and there is another singularity corresponding to the point $t=-\frac{\alpha}{2}$. This singularity is a saddle-node, so, the foliation cannot have a first integral.

Consequently, $2n=k$.
\end{proof}

From now on, we will assume that $k=2n$.
After blowing up $qn$ times the $x$-axis, and $pn$ times the
$y$-axis, the resulting foliation turns out to be, after
eliminating common factors,
$$
\tilde{\Omega}=(2z^2+\alpha zU+2) (pnydx+qnxdy)+xy (2z+\alpha
U)dz.
$$
The only singular points away from the intersections of the
components of the exceptional divisor are those defined by
$(x=0,\ 2z^2+\alpha z+2=0)$. Assume that  $\alpha\neq \pm 4$, so
this equation has two different solutions $\alpha_1$, $\alpha_2$.
Consider the point corresponding $z=\alpha_1$. After a coordinate
change $z-\alpha_1=\tilde{z}$, this point $P_1$ is the origin,
and in the new coordinates the 1-form that generates the
foliation is
$$
\omega= \left( 2\tilde{z} (\tilde{z}+\alpha_1-\alpha_2) +
\alpha (\tilde{z}+\alpha_1) \tilde{U}\right) (pnydx +qnxdy) + xy
(2\tilde{z}-2\alpha_2+\alpha \tilde{U}) d\tilde{z},
$$
where $U=1+\tilde{U}$. The terms of lower order of $\Omega$ are,
in these coordinates,
$$
\Omega_{L,1}= 2xy\tilde{z} (\alpha_1-\alpha_2) \left[ pn
\frac{dx}{x}+qn\frac{dy}{y} +
\frac{\alpha_2}{\alpha_2-\alpha_1}
\frac{d\tilde{z}}{\tilde{z}}\right] .
$$
Similarly, in the coordinates centered around the point $P_2$
corresponding to $z=\alpha_2$, the lower order terms are
$$
\Omega_{L,2}= 2xy\tilde{z} (\alpha_2-\alpha_1) \left[ pn
\frac{dx}{x}+qn\frac{dy}{y} +
\frac{\alpha_1}{\alpha_1-\alpha_2}
\frac{d\tilde{z}}{\tilde{z}}\right] .
$$
In the paper \cite{FM}, the authors studied the generic case
where
$$
\alpha^2\neq \frac{(16+r)^2}{16+2r},
$$
for every $r\in \QQ_{\geq 0}$, where the singularities are
already reduced. We are now interested in the complementary case
$$
\alpha^2= \frac{(16+r)^2}{16+2r},
$$
for some $r\in \QQ_{>0}$. In this case, one of the numbers
$\frac{\alpha_2}{\alpha_2-\alpha_1}$,
$\frac{\alpha_1}{\alpha_1-\alpha_2}$ is a negative rational, and
the other one, a positive rational, this last point being already
reduced. Assume
$\frac{\alpha_2}{\alpha_2-\alpha_1}\in \QQ_{<0}$, and then, $P_2$
is simple resonant. The singularity at $P_1$ may fall in one of
the following two cases:
\be
\item Either it is linearizable, and then, it is dicritical.
\item Or it is of Dulac type.
\ee
In the second case, as we have seen in Section
\ref{singularidadessimples}, the foliation cannot have a
meromorphic first integral. So, we will suppose in the sequel
that we are in  the first case, and then, the foliation is
dicritical. Locally in a neighbourhood of $P_1$, the foliation is
defined by
$$
\Omega_1= xyz\left[ pn\frac{dx}{x}+qn \frac{dy}{y}
-s\frac{d\tilde{z}}{\tilde{z}}\right],
$$
where $s\in \QQ_{>0}$.

In order to linearize $\Omega$ around $P_1$, let us consider the
following vector fields:
\begin{align*}
X_1 & = g(x,y,\tilde{z})\frac{\partial}{\partial
\tilde{z}}-\frac{1}{pn}x\frac{\partial}{\partial x} \\
X_2 & = g(x,y,\tilde{z})\frac{\partial}{\partial
\tilde{z}}-\frac{1}{qn}y \frac{\partial}{\partial y},
\end{align*}
where
$$
g(x,y,\tilde{z})=\frac{2\tilde{z}  (\tilde{z}+\alpha_1-\alpha_2)
+ \alpha (\tilde{z} +\alpha_1) \tilde{U}
}{2 \tilde{z}-2\alpha_2+\alpha \tilde{U}} =
\frac{\alpha_2-\alpha_1}{\alpha_2} z+ h.o.t.
$$
They are commuting and, in the considered case,
$\frac{\alpha_2-\alpha_1}{\alpha_2}\in \QQ_{<0}$. The
linearization is obtained after an analytic transformation in the
$\tilde{z}$ variable: $\tilde{z}=\tilde{z}(x,y,z)$, and then, a
separatrix $\tilde{z}=0$, that perhaps is not the only smooth
separatrix  corresponding to this point. If this is the case, we
choose one of the separatrices. Anyway, there is an infinite
number of separatrices through this point, defined by the
equation
$$
z^{s_1}=Cx^{pns_2}y^{qns_2}, \text{ for }C\in \CC,
$$
where $s=\frac{s_1}{s_2}$. 
Assume that a meromorphic first integral exists, $F(x,y,z)$. In
the coordinates around $P_1$, it is $\tilde{F} (x,y,\tilde{z})$.
Let us follow the same scheme of \cite{MS}: take a transversal
$\Sigma$ to the component of the divisor defined by $y=0$, i.e.,
the special component in the notation of \cite{FM}. For
simplicity, assume that $\Sigma\equiv (x=\tilde{z}=1)$. $\tilde{f}(t)=
\tilde{F}(1,t,1)$ is the restriction of $\tilde{F}$ to $\Sigma$,
defined for small $|t|$. If $t\in \CC$ is small enough, the leaf through
$(1,t,1)$ is
$$
z^2=\frac{1}{t^{qn}} x^{pn}y^{qn},
$$
and so, it approaches any neighbourhood of the origin. Defining
$$
\tilde{f}(t)= F(x,y,(t^{-qn}x^{pn}y^{qn})^{1/s}),
$$
$\tilde{f}(t)$ can be extended to a function on $\bar{\CC}$, necessarily
rational. Let $h_0$, $h_1$ be the generators of the holonomy of
this component, calculated in the transversal $\Sigma$. The fact
that $\tilde{f}(t)$ is the restriction of a meromorphic first integral
implies that $\tilde{f}\circ h_1=h_1$, $\tilde{f}\circ h_2=h_2$, i.e., any
element $h(t)$ of the holonomy group of the special component is
fixed by the rational function $f(t)$.

Conversely, assume that a rational function $f(t)$ exists on a transversal $\Sigma \cong \bar{\CC}$, invariant by the holonomy. Every leaf cuts $\Sigma$ so, following the leaves a first integral may be defined in a neighbourhood of this special component of the exceptional divisor. The other singularities, by \cite{MM}, have holomorphic first integrals, so this meromorphic first integral extends to a neighbourhood of the origin, and defines a meromorphic first integral for the foliation. So, we have proved:
\begin{teorema}
 A meromorphic first integral exists for the quasi-ordinary cuspidal foliation generated by $\Omega= d(z^2+(x^py^q)^k)+\alpha (x^p y^q)^{n} U(x^py^q)dz$ if and only if $k=2n$ and a rational function exists invariant by the projective holonomy of the special component of the divisor that appears after reduction of the singularities as in \cite{FM}.
 \end{teorema}

\section{Dicriticalness} \label{dicriticidad}
Finally, we shall investigate the relation about the dicriticalness  of a nilpotent foliation in $\cen{2}$ with the dicriticalness of its pull-back via a map 
\begin{eqnarray*}
\Phi :(\CC^3,\bo ) & \longrightarrow & (\CC^2,\bo ) \\
(x,y,z) & \longmapsto & (f(x,y), z).
\end{eqnarray*}

So, we will consider cuspidal, nilpotent singularities generated by a 1-form
$$
\Omega= d(z^2+f(x,y)^{n})+\alpha f(x,y)^{p} U(f)dz,
$$
where $U(t)$ is a unit, $U(0)=1$. 

For, let us recall, following \cite{Cano}, that a germ of singular holomorphic foliation $\FF$ in $(\CC^{n},\bo)$, generated by a 1-form $\omega$, is dicritical if there exists an analytic map $\varphi : \cen{2}\rightarrow \cen{n}$, not invariant by $\FF$ (i.e., $\varphi^{\ast}\omega \neq 0$), such that $\varphi^{\ast} \FF$ is the foliation $(dx=0)$, and $\varphi (y=0)$ is invariant by $\FF$.

With this characterization, if $\Omega$ defines a dicritical foliation on $\cen{3}$, also $\omega_{n,p,\alpha,U}$ defines a dicritical foliation in $\cen{2}$. Indeed, if $\varphi : \cen{2}\rightarrow \cen{3}$ exists, in the conditions of the previous characterization, and $\Phi \circ \varphi =\rho$, $\rho^{\ast}\omega_{n,p,\alpha,U}= \varphi^{\ast}\Phi^{\ast} \omega_{n,p,\alpha,U}= \varphi^{\ast}\Omega\neq 0$, $\rho^{\ast}\omega_{n,p,\alpha,U}\equiv (dx=0)$, and $\rho (y=0)=\Phi (\varphi (y=0))$ is invariant.

It is not true, in general, that given a map $\cen{n}\stackrel{\Phi}{\longrightarrow} \cen{m}$, if $\omega$ defines a dicritical foliation in $\cen{m}$, then $\Phi^{\ast} \omega$ also defines a dicritical foliation on $\cen{n}$. In particular, if $m=2$, dicritical foliations are exactly foliations with an infinite number of separatrices, but in higher dimensions, it is not even true that they have a separatrix \cite{J}. So, for a general map $\Phi$, it is not always possible to study the dicriticalness   of the pull-back from the point of view of the existence of infinitely many separatrices. We will see that in our case the situation is different, and much simpler. If $\omega_{n,p,\alpha,U}$ defines a dicritical foliation on $\cen{2}$, then we must have, according to \cite{Meziani,MS}
\be
\item $n=2p$.
\item $\alpha\neq \pm 4$.
\item Under (i) and (ii), after $p$ blowing-ups, the quotient of the eigenvalues at the singular points (different from corners) are, respectively, $\frac{2p (\alpha_{2}-\alpha_1)}{2\alpha_1+\alpha}$, $\frac{2p (\alpha_1-\alpha_2)}{2\alpha_2+\alpha}$, where $\alpha_1$, $\alpha_2$ are the roots of the equation $2y^2+\alpha y +2=0$. It is necessary that one of these numbers is a positive rational. In this case, if it is not an integer or the inverse of an integer, then  it is dicritical. In the remaining situation, either it is dicritical or of Dulac type, depending on the coefficients of $U(x)$.
\ee

Denote, as in the Introduction, $\omega_{p,\alpha, U}=d(z^2+t^{2p})+\alpha
t^{p}U(t)dz$. Let us see that, in our case, if $\omega_{p,\alpha,U}$ is dicritical, then $\Phi^{\ast}\omega_{p,\alpha,U}$  also is. For, we shall use the characterization from \cite{Cano} explicitely. In order to locate a dicritical component in the reduction of singularities of $\omega_{p,\alpha,U}$ (in fact, the only one), we must follow the following steps:
\be
\item Blow-up the origin $p$ times. In an appropriate chart, this is represented by the map $(t,z)\rightarrow (t,t^pz)$.
\item Make a translation $z\mapsto z+\alpha_i$, with previous notations.
\item Blow up the origin a certain number of times, depending on the quotient of the eigenvalues at the singular point. This is done via a map $(t,z)\mapsto (t^{n_1}z^{m_1}, t^{n_2}z^{m_2})$, with $n_1m_2-m_1n_2=1$.
\ee
After these transformations, we arrive at a corner where one of the axis is a dicritical component of the exceptional divisor. The composition of the transformations 1., 2. and 3. is
$$
E(t,z)= \left( t^{n_1}z^{m_1}, t^{pn_1}z^{pm_1} (t^{n_2}z^{m_2}+\alpha_1)\right) .
$$
Assume, for instance, that the horizontal axis $z=0$ is transversal to the foliation $E^{\ast} \FF$ defined by the saturated of $E^{\ast}\omega_{p,\alpha,U}$. In this case, this foliation is generated by $dt+B(t,z)dz$. This foliation can be rectified by a diffeomorphism $S(t,z)=(tS_1 (t,z), z)$, with $S_1 (0,0)\neq 0$. This means that $S^{\ast} (dt+B(t,z)dz)\wedge dx=0$. Assume, first, that we are in the easier case where $f(x,y)=x^{p_1} y^{p_2}$, i.e., the foliation $\Phi^{\ast}\omega_{p,\alpha,U}$ is of quasi-ordinary type. We will call this map $E\circ S$ a \textit{dicriticalness section}, i.e., a map verifying the dicriticalness condition of \cite{Cano} as given above. It is our objective to lift this map to $\cen{3}$ in order to show that $\Omega:= \Phi^{\ast} \omega_{p,\alpha,U}$ is dicritical. The map $E\circ S$ should be lifted to $\cen{3}$ in a map $\sigma=(\sigma_1,\sigma_2,\sigma_3)$ such that $\Phi\circ \sigma = E\circ S$, but this is not always possible, as in particular it would imply that $\sigma_1^{p_1}\sigma_2^{p_2} = t^{n_1} S_1(t,z)^{n_1}z^{m_1}$, which is not possible if $p_1$, $p_2$ are big enough. But, let us observe that, if we define $\varphi  =(\varphi_1,\varphi_2) :\cen{2} \rightarrow \cen{2}$ by $\varphi (t,z) =(t^{r_1}, z^{r_2})$, the composition $E\circ S\circ \varphi $ is again a dicriticalness section for $\omega_{p,\alpha,U}$. Choosing $r_1=p_1$, $r_2=p_2$, this map can be lifted to $\cen{3}$ defining
\begin{align*}
\sigma_1 (t,z) & =t^{n_1}S_1 (t^{p_1},z^{q_1})^{n_1/p_1}, \\ \sigma_2 (t,z) & = z^{m_1}, \\  \sigma_3 (t,z) & = t^{pn_1p_1} S_1 (t^{p_1}, z^{q_1})^{pn_1} z^{pm_1q_1} (t^{n_2p_1} S_1 (t^{p_1}, z^{q_1})^{n_2} z^{m_2q_1} +\alpha_1) .
\end{align*}
So, the foliation $\Phi^{\ast} \omega_{n,p,\alpha}$ is also dicritical.

Let us consider now the general case, i.e., where $\Phi (x,y,z)=(f(x,y),z)$, for a general $f(x,y)$. In this case, consider a local monomialization of $f(x,y)$, i.e., a local map $\rho (x,y): \cen{2}\rightarrow \cen{2}$ such that $f(\rho (x,y))= x^{a}y^{b}V(x,y)$, $V$ being a unit, $a$, $b\neq 0$. As before, modify the dicriticalness section  $E\circ S$ defining $\varphi_1(t,z)=t^a$. This implies that we must choose  $\sigma_1 (t,z)=t^{n_1}$, $\sigma_2 (t,z)=z^{m_1}$. In this case, the equality $f(\rho (\sigma_1,\sigma_2))= \varphi_1 (t,z)^{n_1}S_1 (\varphi_1 (t,z), \varphi_2 (t,z))^{n_1} \varphi_2 (t,z)^{m_1}$ turns out to be 
$$
z^{b} V(t^{n_1}, y^{m_1})^{1/m_1} = S_1 (t^{a}, \varphi_2 (t,z))^{n_1/m_1} \varphi_2 (t,z),
$$
which has a unique solution $\varphi_2 (x,y)$ by Implicit Function Theorem. Now, $\sigma_3 (x,y)$ can be defined accordingly.

This ends the construction. We sum up previous considerations in the following theorem:
\begin{teorema}
The foliation in $\cen{3}$ defined by the 1-form $\Omega= d(z^2+f(x,y)^{2p})+ \alpha f(x,y)^p U(f) dz$ is dicritical if and only if $\omega=d(z^2+t^{2p})+\alpha t^p U(t) dz $ defines a dicritical foliation in $\cen{2}$.
\end{teorema}

\begin{nota}
Previous result also holds in dimension greater than three, as all the arguments are valid there. The differences in the proof for this extended case consist only in technical points, and we shall not reproduce it here.
\end{nota}

Dicritical foliations in dimension two are also characterized for the existence of infinitely many separatrices at the origin. This is not true in higher dimension, where a dicritical foliation may not have any separatrix. Nevertheless, this condition is still sufficient: a foliation with infinitely many separatrices is always dicritical. In the cuspidal case we study in this paper, this situation cannot happen: any nilpotent dicritical foliation generated by a 1-form $\Omega= d(z^2+f(x,y)^{2p})+ \alpha f(x,y)^p U(f) dz$ has infinitely many separatrices. In fact, the foliation in $\cen{2}$ generated by $\omega_{p,\alpha,U}$, in the dicritical case, has separatrices that can be parameterized by  Puiseux series $y=\alpha_1 x^p+\cdots$, series with rational exponents. These series, by composition with $\Phi$, give infinitely many separatrices for $\Omega$.

\end{document}